\documentclass[12pt]{article}
\usepackage{amsxtra,amssymb,amsthm,amsmath,latexsym}

\theoremstyle{plain}

\def\oH{\buildrel\circ\over H}
\def\oH1{\buildrel\circ\over H\kern-.02in{}^1}

\begin{document}


\title{Inequalities for solutions to some nonlinear\\
 equations
   \thanks{key words: inequalities, nonlinear equations, dynamical systems 
method, 
evolution equations
    }
   \thanks{AMS subject classification: 37C35, 37L05, 37N30, 47H15, 58C15, 
65J15, }
}

\author{
A.G. Ramm\\
LMA/CNRS, Marseille 13402, cedex 20, France\\
and  Mathematics Department, 
Kansas State University, \\
 Manhattan, KS 66506-2602, USA\\
ramm@math.ksu.edu\\
}

\date{}

\maketitle\thispagestyle{empty}

\begin{abstract}
Let $F$ be a nonlinear Frechet differentiable map in a
real Hilbert space. Condition sufficient for existence of a solution to 
the equation $F(u)=0$ is given, and a method 
(dynamical systems method, DSM) to calculate the
solution as the limit of the solution to a Cauchy problem is justified
under suitable assumptions.
\end{abstract}


\section{Introduction}
 In this paper  Theorem 1 gives a method for proving existence of the 
solution to nonlinear equation $F(u)=0$ and for 
computing it. Our method 
(dynamical systems method: DSM) consists of solving a suitable Cauchy
problem which has a global solution $u(t)$ such that $y:=u(\infty)$ does 
exist and $F(y)=0$.
 Theorem 1 generalizes Theorem 1 in [1] and its proof is based on the idea
in [1]. The result of Example 1 in Section 3 below was obtained in [5].
The result of Example 5 was essentially obtained in [8], 
where the existence of the solution was assumed, while in Theorem 1 below
the existence of the solution is proved. In Example 6 some of the results
from [4] are obtained.
In Remark 2.4 in [8] some of the 
examples we discuss in Section 3 were mentioned.
In [1] DSM was developed for solving ill-posed operator equations with
not necessarily monotone operators and 
for constructing convergent iterative methods for their solution.
 
Let $F$ be a nonlinear Frechet differentiable map in a real Hilbert 
space. Consider the equation 
$$ F(u)=0.
\eqno{(1.1)}$$
Let $\Phi(t,u)$ be a map, continuous in the norm of $H$ and Lipschitz with 
respect to $u$ in the ball $B:=\{u: ||u-u_0||\leq R, \, u\in H\}$.
Weaker conditions, which guarantee local existence and uniqueness of the 
solution to (1.6) below, would suffice. Assume that:
$$ (F'(u) \Phi (t,u), F(u)) \leq -g_1(t)||F(u)||^2 \quad \forall u\in B,
\eqno{(1.2)}$$
and 
$$ 
||\Phi (t,u)|| \leq g_2(t)||F(u)|| \quad \forall u\in B,
\eqno{(1.3)}$$
where $g_j, j=1,2,$ are positive functions on $R_+:=[0,\infty)$,
$g_2$ is continuous, $g_1\in L_{loc}^1(R_+)$, 
$$ \int_0^\infty g_1dt=+\infty,
\eqno{(1.4)}$$
and
$$ G(t):= g_2(t)\exp (-\int_0^t g_1 ds)\in L^1(R_+).
\eqno{(1.5)}$$
   
{\bf Remark: } Sometimes the assumption (1.3) can be used
in the following modified form:
$$
||\Phi (t,u)|| \leq g_2(t)||F(u)||^b \quad \forall u\in B,
\eqno{(1.3')}
$$
where $b>0$ is a constant. The statements and proofs of Theorems 1-3 
in Sections 1 and 2 can be easily adjusted to this assumption.

Consider the following Cauchy problem:
$$ 
\dot u=\Phi (t,u),   \quad u(0)=u_0,\quad \dot u:=\frac {du}{dt}.
\eqno{(1.6)}$$
Finally, assume that
$$ ||F(u_0)|| \int_0^\infty G(t)dt\leq R.
\eqno{(1.7)}$$
 The above assumptions on $F, \Phi, g_1$ and $g_2$ hold throughout and 
are not repeated in the statement of Theorem 1 below. By global solution 
to (1.6) we mean a solution defined for all $t\geq 0$.

{\bf Theorem 1.} {\it Under the above assumptions problem (1.6)
has a global solution $u(t)$, there exists the 
strong limit $y:=\lim_{t\to \infty} u(t)=u(\infty)$, 
$F(y)=0$, $u(t)\in B$ for all $t\geq 0$,
and the following inequalities hold:

$$ ||u(t)-y||\leq ||F(u_0)||\int_t^\infty G(x)dx,
\eqno{(1.8)}$$

and

$$ ||F(u(t))||\leq ||F(u_0)||\exp(-\int_0^t g_1(x)dx). 
\eqno{(1.9)}$$
}
 Theorem 1 generalizes Theorem 1 in [1].

 In Section 2 proof of Theorem 1 is given, and two other theorems are 
proved, and in Section 3 examples of 
applications are presented.
In Section 4  a linear equation and in Section 5 a nonlinear
equation are discussed.

\section{Proof of Theorem 1 and additional results.}

The assumptions about  $\Phi$ imply local existence and uniqueness of the 
solution $u(t)$ to (1.6). To prove global existence of 
$u$, it is sufficient to prove a uniform with respect to $t$ bound on 
$||u(t)||$. Indeed, if the maximal interval of the existence of $u(t)$
is finite, say $[0,T)$, and  $\Phi(t,u)$ is locally Lipschitz with respect to 
$u$,
then $||u(t)|| \to \infty$ as $t\to T$.
  Let $g(t):=||F(u(t))||$. Since $H$ is real, one uses (1.6)  and (1.2)
to get $g\dot g=(F'(u)\dot u, F)\leq -g_1(t)g^2$, so $\dot g\leq 
-g_1(t)g$, and integrating one gets (1.9), because $g(0)=||F(u_0)||$.
Using (1.3), (1.6) and (1.9), one gets:
$$ ||u(t)-u(s)||\leq g(0)\int_s^t G(x)dx, \quad G(x):=g_2(x)\exp(-\int_0^x 
g_1(z)dz).
\eqno{(2.1)}$$
Because $G\in L^1(R_+)$, it follows from (2.1) that the limit
$y:=\lim_{t\to \infty} u(t)=u(\infty)$ exists, and $y\in B$ by (1.7).
From (1.9) and the continuity of $F$ one gets $F(y)=0$, so $y$ solves 
(1.1).  
Taking $t\to \infty$ and setting $s=t$ in (2.1) yields 
estimate (1.8). The inclusion $u(t)\in B$ for all $t\geq 0$ follows from
(2.1) and (1.7).
Theorem 1 is proved. $\Box$ 

If condition (1.2) is replaced by 
$$ 
(F'\Phi,F)\leq -g_1(t)||F||^a, \quad 0<a<2,
\eqno{(2.2)}$$
then the proof of Theorem 1 yields $g^{1-a}\dot g\leq -g_1(t)$, so
$$ 0\leq g(t) \leq [g^{2-a}(0) -(2-a)\int_0^tg_1(s)ds]^{\frac 1 {2-a}}. 
\eqno{(2.3)}$$
If (1.4) holds, then (2.3) implies $g(t)=0$ for all $t\geq T$, where
$T$ is defined by the equation:
$$ 
g^{2-a}(0) -(2-a)\int_0^T g_1(s)ds=0.
\eqno{(2.4)}$$
Thus $||F(u(t))||=0$ for $t\geq T$. So, by (1.3), $\Phi=0$ for $t\geq T$.
Thus, by (1.6), $u(t)=u(T)$ for $t\geq T$.
Therefore $y:=u(T)$ solves equation (1.1), $F(y)=0$, and
$||u(T)-u(0)||\leq ||F(u_0)||\int_0^Tg_2ds$.
If $||F(u_0)||\int_0^Tg_2ds\leq R$, then $u(t)\in B$  for all $t\geq 0$.
We have proved:

{\bf Theorem 2.} {\it If (1.2) is replaced by (2.2), (1.4) holds,
and $||F(u_0)||\int_0^Tg_2ds\leq R$, where $T$ is defined by (2.4),
then equation (1.1) has a solution in $B=\{u: ||u-u_0||\leq R\}$,
the solution $u(t)$ to (1.6) exists for all $t>0$, $u(t)\in B$, 
$u(t):=y$ for $t\geq T$, and $F(y)=0$, $y\in B$.
} 

If (2.2) holds with $a>2$, and (1.4) holds, then a similar calculation 
yields:
$$ 0\leq g(t) \leq [g^{-(a-2)}(0) +(a-2)\int_0^tg_1(s)ds]^{\frac 1 
{2-a}}:=h(t)\to 0 \quad t\to \infty,
\eqno{(2.5)}$$
because of (1.4).
Assume that
$$ 
\int_0^\infty g_2(s)h(s)ds\leq R.
\eqno{(2.6)}$$
Then (1.3) and (1.6) yield $||u(t)-u(0)||\leq R$,
$||u(t)-u(\infty)||\leq \int_t^\infty g_2(s)h(s)ds\to 0$ as $t\to \infty$.
Therefore an analog of Theorem 1 is obtained:

{\bf Theorem 3.} {\it If (2.2) holds with $a>2$, (1.4) and (2.6) hold,
then the solution $u(t)$ to (1.6) exists for all $t>0$, $u(t)\in B$, there 
exists $u(\infty):=y$ , and $F(y)=0$, $y\in B$.
}

\section{Applications.}

If $g_j=c_j, \, j=1,2,$ and $c_j>0$ are constants, then (1.4) and (1.5) 
hold,
$\int_0^\infty Gdx=c_2/c_1$, so condition (1.7) is:
$$ \frac {c_2}{c_1} ||F(u_0)||\leq R, \quad c_j=c_j(R), \, j=1,2.
\eqno{(3.1)}$$
 Let us give examples of applications of Theorem 1 using
its simplified version with $g_j=c_j >0,\, j=1,2$. 

{\bf Example 1. Continuous Newton-type method.}

Let $\Phi= -[F'(u)]^{-1} F(u)$, and assume 
$$  ||[F'(u)]^{-1}||\leq m_1=m_1(R),\,\,  \forall u\in B.
\eqno{(3.2)}$$
Assumption (3.2) holds in all the examples below. It implies 
"well-posedness" of equation (1.1).
The case of ill-posed equation (1.1), when (3.2) is not valid, was
studied in [2] for linear selfadjoint nonnegative-definite operator,
in [3] for nonlinear monotone hemicontinuous operator, and in [1]
for some other cases.
Then $c_1=1$, $c_2=m_1$,  $\Phi$ is locally Lipschitz if,
for example, one assumes 
$$  ||F''(u)||\leq M_2,\,\,  \forall u\in B,
\eqno{(3.3)}$$
where $M_2=M_2(R)$ is a positive constant, and (3.1) is: 
$$  m_1(R)||F(u_0)||\leq R.
\eqno{(3.4)}$$
In the examples below condition (3.3) is assumed and not repeated.
As we have mentioned in Section 1, this condition can be weakened.

{\bf Conclusion 1:} 
{\it By Theorem 1,  inequality (3.4) implies
existence of a solution $y$ to (1.1) in $B$, global
existence and uniqueness of the solution $u(t)$ to (1.6), convergence  
of $u(t)$ to $y$ as $t \to \infty$,  and the error estimate (1.9).
Condition (3.4) is always satisfied if equation (1.1) has a solution $y$
and if $u_0$ is chosen sufficiently close to $y$.}

{\bf Example 2. Continuous simple iterations method.}

Let $\Phi=-F$, and assume $F'(u)\geq c_1(R)>0$ for all $u\in B$.
Then $c_2=1$, $c_1=c_1(R)$, and (3.1) is:
$$  [c_1(R)]^{-1}||F(u_0)||\leq R.
\eqno{(3.5)}$$
{\it If this inequality holds, then 
Conclusion 1 holds with (3.5) replacing (3.4).}

{\bf Example 3. Continuous gradient method.}

Let $\Phi=-[F']^{*}F$ and assume (3.2).  Then $c_2=M_1(R)$, because
$||[F'(u)]^*||=||F'(u)||\leq M_1(R)$, and 
$(F'\Phi,F)=-||[F'(u)]^*F||^2\leq 
-m_1^{-2}||F||^2$, so $c_1=m_1^{-2}$, where $m_1$ is the constant from 
(3.2).
Here we have used the estimates $||f||= ||A^{-1}Af||\leq ||A^{-1}||||Af||,
||Af||\geq ||A^{-1}||^{-1} ||f||$, with $A:=F'(u)$ and $f=F(u)$.
Estimate (3.1) is: 
$$  M_1m_1^2||F(u_0)||\leq R.
\eqno{(3.6)}$$
{\it If this inequality holds, then 
Conclusion 1 holds with (3.6) replacing (3.4).}

{\bf Example 4. Continuous Gauss-Newton method.}

Let $\Phi=-([F']^{*}F')^{-1}[F']^{*}F$, and assume (3.2). Then
$c_1=1$, $c_2=m_1^2M_1$, and (3.1) is:
$$  M_1m_1^2 ||F(u_0)||\leq R.
\eqno{(3.7)}$$
{\it If this inequality holds, then
Conclusion 1 holds with (3.7) replacing (3.4).}

{\bf Example 5. Continuous modified Newton method.}

$\Phi= -[F'(u_0)]^{-1} F(u)$, and assume 
$ ||[F(u_0)]^{-1}||\leq m_0$. Then $c_2=m_0$. To find $c_1$,
let us note that:
$$(F'\Phi, F)=-||F(u)||^2 -((F'(u)-F'(u_0))[F(u_0)]^{-1}F,F)\leq 0.5 
||F(u)||^2,$$ 
provided that 
$$|((F'(u)-F'(u_0))[F(u_0)]^{-1}F,F)|\leq M_2Rm_0 ||F(u)||^2\leq 0.5 
||F(u)||^2,$$
that is, $M_2Rm_0\leq 0.5$. If  $R=(2M_2m_0)^{-1}$, then the last
inequality becomes an equality.
Choosing such $R$, one has $c_2=m_0$, $c_1=0.5$, and (3.1) is:
$2m_0||F(u_0)||\leq (2M_2m_0)^{-1}$, that is,
$$ 4 m_0^2 M_2||F(u_0)||\leq 1.
\eqno{(3.8)}$$
{\it If this inequality holds, then
Conclusion 1 holds with (3.8) replacing (3.4).}

{\bf Example 6. Descent methods.}

Let $\Phi=-\frac{f}{(f',h)}h$, where $f=f(u(t))$ is a 
differentiable functional  $f: H\to [0,\infty)$, and $h$ is an 
element of $H$.
One has $\dot f=(f',\dot u)=-f$. Thus $f=f_0e^{-t}$, where $f_0:=f(u_0)$.
Assume $||\Phi||\leq c_2|f|^b,\, b>0$. Then $||\dot u||\leq 
c_2|f_0|^be^{-bt}$. Therefore $u(\infty)$ does exist, $f(u(\infty))=0$, 
and $||u(\infty)-u(t)||\leq ce^{-bt}$, $c=const>0$.

If $h=f'$, and $f=||F(u)||^2$, then $f'(u)=2[F']^*(u)F(u)$, 
$\Phi=-\frac{f}{||f'||^2}f'$, and (1.6) is a descent 
method. For this $\Phi$ one has $c_1=\frac 1 2$, and $c_2=\frac {m_1} 
2$, where $m_1$ is defined in (3.2). Condition (3.1) in this example is 
condition (3.4).

{\it If (3.4) holds, then Conclusion 1 holds.}

 We have obtained some results from [4]. Our approach is more general
than the one in [4], since the choices of $f$ and $h$ do not allow one,
e.g., to obtain $\Phi$ used in Example 5. 
  
\section{Remark about linear equations.}

 The following result was proved in [2]: if equation $Ay=f$ in 
a Hilbert space 
has a solution $y$, and $A\geq 0$ is a linear selfadjoint operator, then 
the global solution $u$ to the regularized Cauchy problem 
$$  \quad \dot u=-Au-\alpha(t)u+f, \quad u(0)=u_0,
\eqno{(4.1)}$$
has a limit $\lim_{t\to \infty}u(t):=u(\infty)$, and $A(u(\infty))=f$. 
In [2] $u_0\in H$ is arbitrary, $\alpha>0$ is a continuously 
differentiable,
monotonically decaying to zero as $t\to \infty$, function on $R_+$, 
$\int_0^\infty \alpha dt=+\infty,$ and 
$\alpha^{-2}\dot \alpha \to 0$ as $t\to \infty$. 

Unfortunately the author of [7] was not aware of the paper
[2]. Some of the proofs in [7] are close to these in [2]. 
The author thanks Dr. Ya. Alber for pointing out  references [2] and [3].

If $\alpha >0$, $ \dot \alpha \leq 0,$ and 
$\alpha^{-2}|\dot \alpha|\leq c,$ where $c=const$,  then
$\alpha^{-1}(t)-\alpha^{-1}(0)\leq ct$, so
$\alpha (t)\geq [ct+\alpha^{-1}(0)]^{-1}$ and consequently $\int_0^\infty 
\alpha dt=+\infty$ (cf [7]). Therefore the condition
$\int_0^\infty \alpha dt=+\infty$ in [2] can be dropped.

In this Section we give a new derivation of the result in [2]
under weaker assumptions about $\alpha$, and show that the regularization 
in (4.1) is not necessary.

 First, let us prove that the regularization in (4.1) is not necessary: 
the 
result holds with $\alpha=0$. Below $\to$ denotes strong convergence in 
$H$. 

The solution to (4.1) with $\alpha=0$ is $u(t)=U(t)u_0+\int_0^t 
U(t-s)fds$, where 
$U(t):=\exp (-tA)$. If $E_{\lambda}$ is the resolution 
of the identity of the selfadjoint operator $A$, then
$U(t)u_0=\int_0^{\infty}e^{-t\lambda}dE_{\lambda}u_0\to Pu_0$
as $t\to \infty$, where  $P$ is the
operator of the  orthogonal projection on $N$, and 
 $N$ is the null-space of $A$.
Also $\int_0^t U(t-s)fds= 
\int_0^{\infty} (1-e^{-t\lambda}) dE_{\lambda}y\to y-Py$  
as $t\to \infty$, by the dominated convergence theorem. 
Thus,  $u(\infty)= y-Py+Pu_{0}$ and $A(u(\infty))=f$. $\Box$

Consider now the case $0<\alpha \to 0$ as $t\to \infty$. If
$h(t):=\exp (\int_0^t \alpha(s)ds)$, and $u$ solves (4.1), then 
$$  
u(t)=h^{-1}(t) U(t)u_0+h^{-1}(t)\int_0^{\infty}\exp(-t\lambda) \int_0^t
e^{s\lambda}h(s)ds \lambda dE_{\lambda}y.
\eqno{(4.2)}$$
Using L'Hospital's rule one checks that 
$$\lim_{t \to \infty} \frac {\lambda \int_0^t
e^{s\lambda}h(s)ds}{e^{t\lambda}h(t)}=\lim_{t \to \infty} \frac {\lambda 
e^{t\lambda}h(t)}
{\lambda e^{t\lambda}h(t)+e^{t\lambda}h(t)\alpha(t)}=
1 \quad \forall \lambda >0.
\eqno{(4.3)}$$
From (4.2), (4.3), and the dominated convergence theorem, one gets
$u(\infty)= y-Py$.
The first term on the right-hand side of (4.2) tends to zero as $t\to 
\infty$ (even if $Pu_{0}\neq 0$), if $h(\infty)=\infty$.
To apply the dominated convergence theorem, one checks that
$\frac {\lambda \int_0^t e^{-(t-s)\lambda}h(s)ds}{h(t)}=
\frac {\lambda \int_0^t e^{-s\lambda}h(t-s)ds}{h(t)}
\leq 1$
for all $t>0$ and all $\lambda>0$, where the inequality 
$0<h(t-s)\leq h(t)$, valid for $s\geq 0$, was used. $\Box$

Our derivation uses less
 restrictive assumptions on $\alpha$ than in [2] and [7]: we do not assume
differentiability of $\alpha$, and the property
$\lim_{t\to \infty}\alpha^{-2}\alpha'=0$.
The property $\int_0^\infty \alpha dt=+\infty$, which is equivalent to
$h(\infty)=\infty$, was used above only to prove that
$\lim_{t\to \infty}h^{-1}(t)U(t)u_0= 0$. If 
$\int_0^\infty \alpha 
dt:=q<\infty$, then $u(\infty)=y-Py+e^{-q}Pu_0$, 
and $Au(\infty)=f$, so that the basic conclusion holds
without the assumption $h(\infty)=\infty$.  

Finally, let us prove a typical for ill-posed problems result:
the rate of convergence $u(t)\to y$ can be as slow as one wishes, it is 
not uniform with respect to $f$.
Assume $\alpha =0$, but the proof is 
essentially the same for $0<\alpha \to 0$ as $t\to \infty$.
Assume that $A>0$ is compact, and $A\varphi_j=\lambda_j \varphi_j$,
$(\varphi_j, \varphi_m)=\delta_{jm}$. Then (4.2) with $y=y_m:=\varphi_m$
and $u_0=0$ yields $u(t)=\varphi_m (1-e^{-t\lambda_m})$.
Thus $u(\infty)=y$, but 
for any fixed $T>0$, however large,
one can find $m$ such that $||u(T)-y_m||>0.5$, that is, convergence is 
not 
uniform with respect to $f$.

\section{Remark about nonlinear equations.}

In this Section we give a short and simple proof of the basic
result in [3], and close
a gap in the proof in [3], where it is not
explained why one can apply the L'Hospital rule the second time.

The assumptions in [3] are: {\it the operator $A$ is monotone 
(possibly nonlinear), defined on all of $H$, hemicontinuous, problem (4.1) has 
a unique global solution, 
equation $A(y)=f$ has a solution, $\alpha (t)>0$ decays monotonically to 
zero, $\lim_{t\to \infty} \dot \alpha \alpha^{-2}=0$, and 
$\alpha$ is convex.}

We refer below to these assumptions as A3). If A3) hold, 
the basic result, proved in [3], is the existence of $u(\infty):=y$, 
and the relation $A(y)=f$.
In [3], p. 184, under the additional assumption (1.24) from [3], the 
global 
existence of the solution to (4.1) is proved. Actually,
the assumption about global existence of the solution to (4.1) can be
dropped altogether: in [6], p.99, it is proved that A3) (and even weaker
assumptions)  imply that problem (4.1) has 
a unique global solution.

 Let us give a short
proof of the basic result from [3]. It is well known that A3) 
imply that the problem $A(v_\alpha)+\alpha v_\alpha -f=0,$
for any fixed number $\alpha>0$, has a unique solution,
there exists $\lim_{\alpha 
\to 0} v_\alpha:=y$,  $A(y)=f,$ and $||y||\leq ||z||$, for any
$z\in \{z: A(z)=f\}$. Thus, for any small $\delta>0,$ one can find
$\alpha_\delta$ such that $||v_\alpha-y||<\delta/2$ for all 
$\alpha>\alpha_\delta$, $\lim_{\delta \to 0}\alpha_\delta =0$.
Let $w:=u-v_\alpha$, where $u$ solves (4.1) and $v_{\alpha}$ does 
not depend on $t$. Then $\dot w=-[A(u)-A(v_\alpha) +\alpha(t)(u-v_\alpha) 
+(\alpha(t)-\alpha) v_\alpha]$. Multiply this by $w$, use the monotonicity 
of $A$, and let $||w||:=g$.
Then $g\dot g\leq -\alpha(t)g^2 +c|\alpha(t)-\alpha|g$, $c=||y||$.
Indeed, multiply $A(v_\alpha)+\alpha v_\alpha-A(y)=0$ by $v_\alpha-y$
and use monotonicity of $A$ to get $\alpha (v_\alpha,v_\alpha-y)\leq 0$. 
Thus
$||v_\alpha||\leq ||y||$, so $c=||y||$. 

Since $\alpha(t)$ is convex,
one has $|\alpha(t)-\alpha|\leq |\dot \alpha(t)|(t_\alpha -t)$,
where $t_\alpha\geq t$ is defined by the equation 
$\alpha=\alpha (t_\alpha)$, $\lim_{\alpha \to 0}t_\alpha=\infty$. 
Thus,
  $g\dot g\leq -\alpha(t)g^2 +c|\dot \alpha(t)|(t_\alpha -t)g$, and,
taking $u(0)=v_\alpha$, one gets
$$
g(t_\alpha)\leq c e^{-\int_0^{t_\alpha}\alpha (x)dx}
\int_0^{t_\alpha} e^{\int_0^s\alpha (x)dx}|\dot \alpha(s)|(t_\alpha -s)ds.
\eqno{(5.1)}
$$
We prove below that
$$  
\lim_{t\to \infty} \alpha (t) e^{\int_0^t\alpha(s)ds}=\infty.
\eqno{(5.2)}
$$
This allows one to apply twice L'Hospital's rule to the right-hand side of 
(5.1), and to get: 
$\lim_{\alpha \to 0} g(t_\alpha)=\lim_{t_{\alpha} \to \infty} \frac {\dot 
\alpha 
(t_\alpha)}{\dot \alpha (t_\alpha)
+\alpha^2(t_\alpha)}=0$. 
Now, $||u(t_\alpha)-y||\leq ||u(t_\alpha)-v_\alpha||+||v_\alpha-y||$,  
$||v_\alpha-y||\leq \delta/2$, and 
$||u(t_\alpha)-v_\alpha||\leq \delta/2$, 
for sufficiently large $t_\alpha$. Since $\delta>0$ is 
arbitrarily small, it follows that $\lim_{t \to \infty}
||u(t)-y||=0$. 

Let us prove (5.2). From our assumptions about $\alpha$, it follows that 
for all sufficiently large $t$, one has $-\dot \alpha \alpha^{-2}\leq c$,
where $0<c<1$, so $\alpha (t)\geq (c_1+t)^{-1} b$, where $b:=c^{-1}>1$, 
$c_1>0$ is a constant,
and $e^{\int_0^t\alpha(s)ds}\geq (c_1+t)^b$. Thus, (5.2) holds.
$\Box$

If one assumes additionally that $A$ is Frechet differentiable, then the 
proof is shorter. Namely, let 
$h(t):=||A(u(t))+\alpha(t)u(t) -f||:=||\psi||.$ Then
$h\dot h=-((A'(u(t))+\alpha(t))\psi, \psi)\leq -\alpha(t)h^2$, because
$A'\geq 0$ due to the monotonicity of $A$. Thus $h(t)\leq  \phi (t)$,
where $\phi (t):=h(0) e^{-\int_0^t\alpha ds}$. As we proved in Section 4, 
the 
assumptions on $\alpha (t)$ imply $\alpha (t)\geq (c_1t+c_2)^{-1}$,
where $c_1$ and $c_2$ are positive constants, and $c_1$ can be chosen
so that $0<c_1<1$, due to the assumption 
$\lim_{t\to \infty} \dot \alpha \alpha^{-2}=0$. Therefore
$\int_0^\infty \phi(t) dt<\infty$. From (4.1) one gets: $||\dot u||\leq 
\phi(t)$. Because $\int_0^\infty \phi(t) dt<\infty$, it follows that
$u(\infty):=y$ exists, and $||u(\infty)-u(t)||\leq \int_t^\infty \phi 
(s)ds$. Finally, $A(y)=f$ because $h(\infty)=0=||A(y)-f||$.
Any choice of $\alpha$, for which $\int_0^\infty \phi(t) dt<\infty$, 
is sufficient for the above argument.  $\Box$

\end{document}